\documentclass{amsart}

\usepackage{url,graphicx}
\usepackage[dvips]{epsfig}
\usepackage{latexsym,color}
\usepackage{amscd,amssymb,verbatim}

\newcommand{\be}{\begin{enumerate}}
\newcommand{\ee}{\end{enumerate}}

\newcommand{\pr}{\nin{\bf Proof. }}
\newcommand{\nin}{\noindent}

\newcommand{\br}{{\mathbb R}}

\newcommand{\ca}{{\mathcal A}}

\newcommand{\cf}{{\mathcal F}}

\newcommand{\co}{{\mathcal O}}
\newcommand{\catr}{{\mathcal R}}
\newcommand{\cs}{{\mathcal S}}

\newcommand{\ctop}{{\bf Top}}

\newcommand{\dar}{\downarrow}
\newcommand{\Dar}{\Downarrow}

\newcommand{\id}{\text{id}}

\newcommand{\hra}{\hookrightarrow}

\newcommand{\im}{\text{Im}\,}

\newcommand{\lcm}{\text{lcm}\,}
\newcommand{\lra}{\longrightarrow}

\newcommand{\minus}{\text{Minus\,}}

\newcommand{\plus}{\text{Plus\,}}

\newcommand{\ra}{\rightarrow}
\newcommand{\rsa}{\rightsquigarrow}

\newcommand{\slx}{\Sigma_\lambda^X}

\newcommand{\lspan}{\text{span}\,}
\newcommand{\susp}{\text{susp}\,}

\newcommand{\thra}{\twoheadrightarrow}
\newcommand{\ti}{\tilde}
\newcommand{\tspan}{{\text{span}\,}}

\newcommand{\un}{\text{un}\,}

\newtheorem{thm}{Theorem}[section]
\newtheorem{df}[thm]{Definition}
\newtheorem{lm}[thm]{Lemma}
\newtheorem{crl}[thm]{Corollary}
\newtheorem{prop}[thm]{Proposition}

\newtheorem{rem}[thm]{Remark}

\newtheorem{exam}[thm]{Example}
\numberwithin{equation}{section}

\begin{document}

        \title{Resonance Category}
              \author{Dmitry N. Kozlov}
                   \date{\noindent\today \\[0.05cm] 
\mbox{ }  \hskip6pt  Mathematics Subject Classification (2000):
   Primary 32S20, Secondary 18B30, 32S60, 58K15. \\[0.05cm]
\mbox{ }  \hskip6pt Keywords: spaces of polynomials, symmetric 
smash products, stratifications, resonances,  \\[0.05cm] 
\mbox{ }  \hskip6pt partitions.\\[0.05cm] 
\mbox{ }  \hskip6pt  This research was supported by the Research Grant  
of the Swiss National Science Foundation.
}

\address{ Department of Mathematics, Royal Institute of Technology,
S-100 44, Stockholm, Sweden.}

\email{kozlov@math.kth.se.} 

\begin{abstract} 
  The main purpose of this paper is to introduce a~new category, 
which we call a~resonance category, whose combinatorics reflect that 
of canonical stratifications of $n$-fold symmetric smash products.
The study of the stratifications can then be abstracted
to the study of functors satisfying certain sets of axioms, 
which we name resonance functors.

One frequently studied stratification is that of the set of all
polynomials of degree $n$, defined by fixing the allowed
multiplicities of roots. We apply our abstract combinatorial
framework, in particular, the notion of direct product of relative
resonances, to study the Arnold problem of computing the
algebro-topological invariants of these strata.
\end{abstract}

\maketitle

                 \section{Introduction}
              
 Complicated combinatorial problems often arise when one studies the
homological properties of strata in some topological space with
a~given natural stratification. In this paper, we study the symmetric
smash products stratified by point multiplicities.

More specifically, let $X$ be a~pointed topological space (we refer to
the base point as a~point at infinity), and denote
$$X^{(n)}=\overbrace{X\wedge X\wedge\dots\wedge X}^n/\cs_n,$$ where
$\wedge$ is the smash product of pointed spaces.  In other words,
$X^{(n)}$ is the set of all unordered collections of $n$ points on $X$
with the collections having at least one of the points at infinity
identified, to form a~new infinity point. $X^{(n)}$ is naturally
stratified by point coincidences, and the strata are indexed by the
number partitions of $n$. Note that we consider the closed strata, so,
for example, the stratum indexed $(\underbrace{1,1,\dots,1}_n)$ is the
whole space $X^{(n)}$.

The main open stratum, that is the complement of the closed stratum
$(2,1,\dots,1)$, is a~frequently studied object. It was suggested by
Arnold in a~much more general context, see for example~\cite{Ar1},
that in situations of this kind one should study the problem for all
closed strata. The main argument in support of this point of view is
that there is usually no natural stratification on the main open
stratum, while there is one on its complement, also known as
discriminant.  Having a~natural stratification allows one to apply
such computational techniques, as spectral sequences, in a~canonical
way. Once some information has been obtained about the closed strata,
one can try to find out something about the open stratum by means of
some kind of duality.

If one specifies $X=S^1$, resp.~$X=S^2$, one obtains as strata the
spaces of all monic real hyperbolic, resp.~monic complex, polynomials
of degree~n with specified root multiplicities. These spaces naturally
appear in singularity theory, \cite{AGV}. Homological invariants of
several of these strata were in particular computed by Arnold,
Shapiro, Sundaram, Welker, Vassiliev, and the author,
see~\cite{Ar1,Ko1,Ko2,SW,V2}. These are the special cases which have
inspired this general study.
  
Here, we take a different, more abstract look at this set of
problems. More specifically, the idea is to introduce a~new canonical
combinatorial object, independent of topology of particular $X$, where
the combinatorial aspects of these stratifications would be fully
reflected. This object is a~certain category, which we name the {\bf
resonance category}. It was suggested to the author by
B.~Shapiro,~\cite{Sh}, to use the term {\it resonance} as a~generic 
reference to a~certain type of linear relations among parts of
a~number partition.

Having this canonically defined category at hand, one then can, for
each specific topological space $X$, view the natural stratification
of $X^{(n)}$ as a~certain functor from the resonance category to
$\ctop^*$. These functors satisfy a~system of axioms, which we take as
a~definition of {\bf resonance functors}.  The combinatorial
structures in the resonance category will then project to the
corresponding structures in each specific $X^{(n)}$. This opens the
door to develop the general combinatorial theory of the resonance
category, and then prove facts valid for all resonance functors
satisfying some further conditions, such as for example acyclicity of
certain spaces.

The main combinatorial structure inside the resonance category, which
we study, is that of {\it relative resonances} and their direct
products. Intuitively, a~relative resonance encodes the combinatorial
type of a~stratum with a~union of some substrata shrunk to form the
new infinity point. These spaces appear naturally if we are trying to
compute the homology groups of our strata by means of long exact
sequences, or, more generally, spectral sequences.

Our idea is that the combinatorial knowledge of which relative
resonances are {\it reducible} (that is, are direct products of other
relative resonances), serves as a~guidance for which long exact
sequences one is to consider for the actual homology computations.
This way, the Arnold problem of computing the algebraic invariants of
the strata, splits into two parts: the combinatorial one, embodied by
various structures in the resonance category, such as the relative
resonances, and the topological one, reflecting the specific
properties of~$X$.

Our notions of sequential and strongly sequential resonances are
intended to capture the combinatorial structure of those resonances,
which are particularly compatible with the spectral sequence
computations. This, in turn, leads to the natural notion of complexity
of resonances.

As mentioned above, to illustrate a~possible appearance of this
abstract framework we choose to use a~class of topological spaces
which come in particular from the singularity theory, and whose
topological properties have been studied: spaces of polynomials (real
or complex) with prescribed root multiplicities. In particular, in
case of strata $(k^m,1^t)$, which were studied in \cite{Ar1,Ko1}
for the complex case, and in \cite{Ko2,SW} for the real case, we
demonstrate how the inherent combinatorial structure of the resonance
category makes this particular resonance especially reducible.

\vskip4pt
\nin
The paper is organized as follows: 

\nin
{\bf Section 2.} We introduce the notion of resonance category, 
and describe the structure of its set of morphisms.

\nin
{\bf Section 3.} We introduce the notions of relative resonances,
direct products of relative resonances, and resonance functors.

\nin {\bf Section 4.} We formulate the problem of Arnold and Shapiro
which motivated this research as that concerning a~specific resonance
functor. Then, we analyze the combinatorial structure of resonances
$(a^k,b^l)$, which leads to the complete determination of the homotopy
types of the corresponding strata for $X=S^1$.

\nin {\bf Section 5.} We analyze the combinatorial structure of the
sequential and strongly sequential resonances. For $X=S^1$, this leads
to the complete computation of homotopy types of the strata
corresponding to resonances $(a^k,b^l,1^m)$, such that $a-bl\leq
m$. Next, we consider division chain resonances, which constitute
a~vast generalization of the case $(a^k,1^l)$. We prove, that in this
case the strata always have a homotopy type of a~bouquet of
spheres. We describe a~combinatorial model to enumerate these spheres
as paths in a certain weighted directed graph, with dimensions of the
spheres being given by the total weights of the paths.

\nin
{\bf Section 6.} We introduce the notion of a complexity of a~resonance
and give a~series of examples of resonances having arbitrarily high
complexity. 

\vskip4pt \nin 
{\bf Acknowledgments.}  I would like to thank Peter Mani-Levitska and
Eva-Maria Feichtner for the helpful discussions during the preparation
of this paper.  I am grateful to the anonymous referee, whose comments
helped to make this paper more transparent.  I also express my
gratitude to the Swiss National Science Foundation for supporting this
research.

         \section{Resonance Category}

\subsection{Resonances and their symbolic notation.} $\,$

    For every positive integer $n$, let $\{-1,0,1\}^n$ denote the set of 
all points in $\br^n$ with coordinates in the set $\{-1,0,1\}$. 
We say that a subset $S\subseteq\{-1,0,1\}^n$ is {\bf span-closed} if 
$\tspan(S)\cap\{-1,0,1\}^n=S$, where $\tspan(S)$ is the linear subspace 
spanned by the origin and points in $S$. Of course the origin lies 
in every span-closed set. For $x=(x_1,\dots,x_n)\in\{-1,0,1\}^n$,
we use the notations $\plus(x)=\{i\in[n]\,|\,x_i=1\}$ and 
$\minus(x)=\{i\in[n]\,|\,x_i=-1\}$. 

\begin{df} $\,$
  
  \nin (1) A subset $S\subseteq\{-1,0,1\}^n$ is called an~{\bf
    $n$-cut} if it is span-closed and for every $x\in
  S\setminus\{$origin$\}$ we have $\plus(x)\neq\emptyset$ and
  $\minus(x)\neq\emptyset$.  We denote the set of all
  $n$-cuts by $\catr_n$.

\nin
(2) $\cs_n$ acts on $\{-1,0,1\}^n$ by permuting coordinates, which in turn
induces $\cs_n$-action on $\catr_n$. The {\bf $n$-resonances} are defined
to be the orbits of the latter $\cs_n$-action. We let $[S]$ denote the
$n$-resonance represented by the $n$-cut $S$.
\end{df}

The resonance consisting of origin only is called {\it trivial}.

\begin{exam}\label{exsmres}
{\bf $n$-resonances for small values of $n$.}

(1) There are no nontrivial 1-resonances.

(2) There is one nontrivial 2-resonance: $[\{(0,0),(1,-1),(-1,1)\}]$.

(3) There are four nontrivial 3-resonances:
  $$[\{(0,0,0),(1,-1,0),(-1,1,0)\}],$$
  $$[\{(0,0,0),(1,-1,0),(-1,1,0),(1,0,-1),(-1,0,1),(0,1,-1),(0,-1,1)\}],$$
  $$[\{(0,0,0),(1,-1,-1),(-1,1,1)\}],$$
  $$[\{(0,0,0),(1,-1,-1),(-1,1,1),(0,1,-1),(0,-1,1)\}].$$ 

(4) Here is an example of a nontrivial 6-resonance:
$$[\{(0,0,0,0,0,0),\pm(1,1,0,-1,-1,0),\pm(0,1,1,0,-1,-1),\pm(1,0,-1,-1,0,1)\}].$$
\end{exam}

\nin
{\bf Symbolic notation.} To describe an $n$-resonance, rather than to list
all of the elements of one of its representatives, it is more convenient
to use the following symbolic notation: we write a sequence of $n$ linear
expressions in some number (between 1 and $n$) of parameters, the order
in which the expressions are written is inessential.

Here is how to get from such a~symbolic expression to the
$n$-resonance: choose an order on the $n$ linear expressions and
observe that now they parameterize some linear subspace of ${\mathbb R}^n$, 
which we denote by $A$. The $n$-resonance is now the orbit of
$A^{\perp}\cap\{-1,0,1\}^n$.

Reversely, to go from an $n$-resonance to a~symbolic expression: choose
a~representative $n$-cut $S$, the symbolic expression can now 
be obtained as a~linear parameterization of $\tspan(S)^\perp$.

For example the 6 nontrivial resonances listed in the
Example~\ref{exsmres} are (in the same order):
$$(a,a),\,(a,a,b),\,(a,a,a),\,(a+b,a,b),\,(2a,a,a),\,
(a+b,b+c,a+d,b+d,c+d,2d).$$

\subsection{Acting on cuts with ordered set partitions.} $\,$
 
We say that $\pi$ is an {\bf ordered set partition} of $[n]$ with $m$
parts (sometimes called {\it blocks}) when $\pi=(\pi_1,\dots,\pi_m)$,
$\pi_i\neq\emptyset$, $[n]=\cup_{i=1}^m\pi_i$, and
$\pi_i\cap\pi_j=\emptyset$, for $i\neq j$.  If the order of the parts
is not specified, then $\pi$ is just called a~{\bf set partition}. We
denote the set of all partitions, resp.\ ordered partitions, of a set
$A$ by $P(A)$, resp.\ $OP(A)$. For $P([n])$, resp.\ $OP([n])$, we use
the shorthand notations $P(n)$, resp.\ $OP(n)$. Furthermore, for every
set $A$, we let $\un:OP(A)\ra P(A)$ be the map which takes the ordered
partition to the associated unordered partition.

\begin{df} 
  Given $\pi=(\pi_1,\dots,\pi_k)$ an ordered set partition of $[m]$
  with $k$ parts, and $\nu=(\nu_1,\dots,\nu_m)$ an ordered set
  partition of $[n]$ with $m$ parts, their {\bf composition}
  $\pi\circ\nu$ is an ordered set partition of $[n]$ with $k$ parts,
  defined by $\pi\circ\nu=(\mu_1,\dots,\mu_k)$,
  $\mu_i=\cup_{j\in\pi_i}\nu_j$, for $i=1,\dots,k$.
  
  \nin Analogously, we can define $\pi\circ\nu$ for an~ordered set
  partition $\nu$ and a~set partition $\pi$, in which case
  $\pi\circ\nu$ is a~set partition without any specified order on the
  blocks.
\end{df}

In particular, when $m=n$, and $|\pi_i|=1$, for $i=1,\dots,n$, we can
identify $\pi=(\pi_1,\dots,\pi_n)$ with the corresponding permutation
of $[n]$. The composition of two such ordered set partitions
corresponds to the multiplication of corresponding permutations, and
we denote the ordered set partition $(\{1\},\dots,\{n\})$ by $\id_n$,
or just~$\id$.

\begin{df}
  For $A\subseteq B$, let $p_{B,A}:P(B)\ra P(A)$ denote map induced by
  the restriction from $B$ to $A$. For two disjoint set $A$ and $B$,
  and $\Pi\subseteq P(A)$, $\Lambda\subseteq P(B)$, we define
  $\Pi\times\Lambda=\{\pi\in P(A\cup B)\,|\,p_{A\cup B,A}(\pi)\in\Pi,
  p_{A\cup B,B}(\pi)\in\Lambda\}$.

\end{df}

The following definition provides the combinatorial constructions
necessary to describe the morphisms of the resonance category, as well
as to define the relative resonances.

\begin{df} \label{df2.4}
Assume $S$ is an~$n$-cut. 
For an~ordered set partition of $[n]$, denoted
$\pi=(\pi_1,\dots,\pi_m)$, we define $\pi S\in\catr_m$ to be the set
of all $m$-tuples $(t_1,\dots,t_m)\in\{-1,0,1\}^m$, for which there
exists $(s_1,\dots,s_n)\in S$, such that for all $j\in[m]$, and
$i\in\pi_j$, we have $s_i=t_j$.
\end{df}

\nin Clearly $\id S=S$, and one can see that $(\pi\circ\nu)S=\pi(\nu S)$.

\vskip3pt

\nin {\bf Verification of $(\pi\circ\nu) S=\pi(\nu S)$.}

\nin By definition we have
$$ (\pi\circ\nu)S=\{(t_1,\dots,t_k)\,|\,\exists(s_1,\dots,s_n)\in S
\text{ s.t. }\forall j\in[k],i\in\mu_j:s_i=t_j\},$$
$$\nu S=\{(x_1,\dots,x_m)\,|\,\exists(s_1,\dots,s_n)\in S
\text{ s.t. }\forall q\in[m],i\in\nu_q:s_i=x_q\},$$
$$\pi(\nu S)=\{(t_1,\dots,t_k)\,|\,\exists(x_1,\dots,x_m)\in\nu S
\text{ s.t. }\forall j\in[k],q\in\pi_j,i\in\nu_q:s_i=t_j\}.$$
The identity $(\pi\circ\nu) S=\pi(\nu S)$ follows now from 
the equality $\mu_j=\cup_{q\in\pi_j}\nu_q$.

\vskip3pt

There are many different ways to formulate the Definition~\ref{df2.4}.
We chose the ad hoc combinatorial language, but it is also possible to
put it in the linear-algebraic terms. An ordered set partition of
$[n]$, $\pi=(\pi_1,\dots,\pi_m)$, defines an inclusion map
$\phi:\br^m\ra\br^n$ by $\phi(e_i)= \sum_{j\in\pi_i}\tilde e_j$, where
$\{e_1,\dots,e_m\}$, resp.~$\{\tilde e_1,\dots,\tilde e_n\}$, is the
standard orthonormal basis of $\br^m$, resp.~$\br^n$. Given
$S\in\catr_n$, $\pi S$ can then be defined as
$\phi^{-1}(\text{Im}\,\phi\cap S)$.

\subsection{The definition of the resonance category and 
the terminology for its morphisms.} $\,$

\begin{df}
  The {\bf resonance category}, denoted $\catr$, is defined as follows:
  
  \nin (1) The set of objects is the set of all $n$-cuts, for
  all positive integers $n$, $\co(\catr)=\cup_{n=1}^\infty \catr_n$.
  
  \nin (2) The set of morphisms is indexed by triples $(S,T,\pi)$,
  where $S\in\catr_m$, $T\in\catr_n$, and $\pi$ is an ordered set
  partition of $[n]$ with $m$ parts, such that $S\subseteq \pi T$.
  For the reasons which will become clear later we denote the morphism
  indexed with $(S,T,\pi)$ by $S\thra\pi T\stackrel{\pi}{\hra} T$.

\nin
As the notation suggests, the initial object of the morphism 
$S\thra\pi T\stackrel{\pi}{\hra} T$ is $S$ and terminal object is $T$.
The composition rule is defined by
$$(S\thra\pi T\stackrel{\pi}{\hra}T)\circ(T\thra\nu
Q\stackrel{\nu}{\hra}Q)= S\thra\pi\nu Q\stackrel{\pi\nu}{\hra} Q,$$
where $S\in\catr_k$, $T\in\catr_m$, $Q\in\catr_n$, $\pi$ is an ordered
set partition of $[m]$ with $k$ parts, and $\nu$ is an ordered set
partition of $[n]$ with $m$ parts.
\end{df}

An~alert reader will notice that the resonances themselves did not
appear explicitly in the definition of the resonance category. In
fact, it is not difficult to notice that resonances are isomorphism
classes of objects of $\catr$. Let us now look at the set of morphisms
of $\catr$ in some more detail.

\nin
(1) For $S\in\catr_n$, the identity morphism of $S$ is
 $S\thra S\stackrel{\id}{\hra} S$. 

\nin
(2) Let us introduce short hand notations: $S\thra T$ for
$S\thra T\stackrel{\id}{\hra} T$, and $\pi T\stackrel{\pi}{\hra}T$
for $\pi T\thra \pi T\stackrel{\pi}{\hra}T$. Then we have
$$S\thra\pi T\stackrel{\pi}{\hra} T=
(S\thra\pi T)\circ(\pi T\stackrel{\pi}{\hra}T).$$
Note also that $S\thra S=S\stackrel{\id}{\hra}S$.

\nin
(3) The associativity of the composition rule can be derived from 
the commutation relation 
$$(\pi S\stackrel{\pi}{\hra}S)\circ(S\thra T)=
(\pi S\thra\pi T)\circ(\pi T\stackrel{\pi}{\hra} T)$$
as follows:
\begin{multline} \notag
(S\thra\pi T\hra T)\circ(T\thra\nu Q\hra Q)\circ(Q\thra\rho X\hra X)=\\
(S\thra\pi T)\circ(\pi T\hra T)\circ(T\thra\nu Q)\circ(\nu Q\hra Q)
\circ(Q\thra\rho X)\circ(\rho X\hra X)=\\
(S\thra\pi T)\circ(\pi T\thra\pi\nu Q)\circ(\pi\nu Q\thra\pi\nu\rho X)\circ \\
(\pi\nu\rho X\hra \nu\rho X)\circ(\nu\rho X\hra \rho X)\circ(\rho X\hra X).
\end{multline}

\nin
(4) We shall use the following names: morphisms $S\thra T$ are called
{\bf gluings} (or $n$-gluings, if it is specified that $S,T\in\catr_n$); 
morphisms $\pi T\stackrel{\pi}{\hra}T$ are called {\bf inclusions}
(or $(n,m)$-inclusions, if it is specified that $T\in\catr_n$, 
$\pi T\in\catr_m$), the inclusions are called {\bf symmetries} if $\pi$ 
is a~permutation. As observed above, the symmetries are the only 
isomorphisms in $\catr$. Here are two examples of inclusions:
$$\{(0,0),(1,-1),(-1,1)\}\stackrel{(\{1\},\{2,3\})}{\hookrightarrow}
\{(0,0,0),(1,-1,-1),(-1,1,1)\},$$
$$\{(0,0),(1,-1),(-1,1)\}\stackrel{(\{1\},\{2,3\})}{\hookrightarrow}
\{(0,0,0),\pm(1,-1,-1),\pm(0,1,-1)\}.$$

\section{Relative resonances, direct products, and resonance functors.}

\subsection{Relative Resonances.} $\,$

Let $A(n)$ denote the set of all collections of non-empty multisubsets
of~$[n]$, and let $P(n)\subseteq A(n)$ be the set of all partitions
of~$[n]$. For every $S\in\catr_n$ let us define a~closure operation on
$A(n)$, resp.\ on $P(n)$.

\begin{df}\label{df3.1}
Let $\ca\in A(n)$. We define $\ca\Dar S\subseteq A(n)$
to be the minimal set satisfying the following conditions:
\begin{enumerate}
\item $\ca\in\ca\Dar S$;
\item if $\{B_1,B_2,\dots,B_m\}\in\ca\Dar S$, then 
$\{B_1\cup B_2,B_3,\dots,B_m\}\in\ca\Dar S$;
\item if $\{B_1,B_2,\dots,B_m\}\in\ca\Dar S$, and there exists 
$x\in S$, such that $\plus(x)\subseteq B_1$, then
$\{(B_1\setminus\plus(x))\cup\minus(x),B_2,\dots,B_m\}\in\ca\Dar S$.
\end{enumerate}

For $\pi\in P(n)$, we define $\pi\dar S\subseteq P(n)$ as 
$\pi\dar S=(\pi\Dar S)\cap P(n)$. For a set $\Pi\subseteq P(n)$
we define $\Pi\dar S=\cup_{\pi\in\Pi}\pi\dar S$. We say that $\Pi$
is {\bf $S$-closed} if $\Pi\dar S=\Pi$.
\end{df}

The idea behind this definition comes from the context of the standard
stratification of the $n$-fold symmetric product. Given a~stratum $X$
indexed by a number partition of $n$ with $m$ parts, let us fix some
order on the parts. A~substratum $Y$ is obtained by choosing some
partition $\pi$ of $[m]$ and summing the numbers within the blocks of
$\pi$. Since the order of the parts of the number partition indexing
$X$ is fixed, $X$ gives rise to a~unique $m$-cut $S$. The set
$\pi\dar S$ describes all partitions $\nu$ of $[m]$ such that if the
numbers within the blocks of $\nu$ are summed then the obtained
stratum $Z$ satisfies $Z\subseteq Y$. In particular, if $Y$ is shrunk
to a~point, then so is $Z$. The two following examples illustrate how
the different parts of the Definition~\ref{df3.1} might be needed.

\vskip3pt

\noindent
{\bf Example 1.} {\it The equivalences of type (2) from
  the~Definition~\ref{df3.1} are needed.} Let the stratum $X$ be
indexed by $(3,2,1,1,1)$ (fix this order of the parts), and let
$\pi=\{1\}\{23\}\{4\}\{5\}$. Then, the stratum~$Y$ is indexed by
$(3,3,1,1)$. Clearly, the stratum~$Z$, which is indexed by $(3,3,2)$,
lies inside $Y$, hence $\{1\}\{2\}\{345\}\in\pi\dar S$, where $S$ is
the cut corresponding to $(3,2,1,1,1)$. However, if one
starts from the partition $\pi$ and uses equivalences of type (3) from
the Definition~\ref{df3.1}, the only other partitions one can obtain
are $\{1\}\{24\}\{3\}\{5\}$, and $\{1\}\{25\}\{3\}\{4\}$. None of them
refines $\{1\}\{2\}\{345\}$, hence it would not be enough in the
Definition~\ref{df3.1} to just take the partitions which can be
obtained via the equivalences of type (3) and then take $\pi\dar S$ to
be the set of all the partitions which are refined by these.

\vskip3pt

\noindent
{\bf Example 2.} {\it It is necessary to view the equivalence relation
  on the larger set $A(n)$.} This time, let the stratum $X$ be indexed
by $(a+b,b+c,a+d,b+d,c+d,2d)$ (fix this order of the parts, and assume
as usual that there are no linear relations on the parts other than
those induced by the algebraic identities on the variables $a$, $b$,
$c$, and $d$). Furthermore, let $\pi=\{16\}\{23\}\{45\}$. Then the
stratum $Y$ is indexed by $(a+b+2d,a+b+c+d,b+c+2d)$. Clearly, we have
$\{34\}\{15\}\{26\}\in\pi\dar S$, where $S$ is the cut
corresponding to $(a+b,b+c,a+d,b+d,c+d,2d)$. 

\noindent
A natural idea for the
Definition~\ref{df3.1} could have been to define the equivalence
relation directly on the set $P(n)$ and use ``swaps'' instead of the
equivalences of type~(3), i.e., to replace the condition (3) by:

\begin{quote} {\it
  if $\{B_1,B_2,\dots,B_m\}\in\ca\dar S$, and there exists $x\in S$,
  such that $\plus(x)\subseteq B_1$, and $\minus(x)\subseteq B_2$,
  then $\{(B_1\setminus\plus(x))\cup\minus(x),
  (B_2\setminus\minus(x))\cup\plus(x),B_3,\dots,B_m\}\in\ca\dar S$. }
\end{quote}

\noindent
However, this would not have been sufficient as this example shows,
since no swaps would be possible on $\pi=\{16\}\{23\}\{45\}$.

\begin{df}
  Let $S$ be an~$n$-cut, $\Pi\subseteq P(n)$ an~$S$-closed
  set of partitions. We define
  $$S\setminus\Pi=S\setminus\{s\in S\,|\,(\minus(s),\plus(s),
\text{singletons})\in\Pi\},$$
where $(\minus(s),\plus(s),\text{singletons})$ is the partition
which has only two nonsingleton blocks: $\minus(s)$ and $\plus(s)$.  
\end{df}

In the next definition we give a~combinatorial analog of viewing
a~stratum relative to a~substratum.

\begin{df} \label{dfrel} $\,$
  
  \nin (1) A {\bf relative $n$-cut} is a~pair $(S,\Pi)$,
where $S\subseteq\{-1,0,1\}^n$, $\Pi\subseteq P(n)$, such that the 
following two conditions are satisfied:
\begin{itemize}
\item $(\lspan S)\setminus\Pi=S$;
\item $\Pi$ is $(\lspan S)$-closed.
\end{itemize}
  
\nin (2) The permutation $\cs_n$-action on $\{-1,0,1\}^n$ induces
an~$\cs_n$-action on the relative $n$-cuts by
$(S,\Pi)\stackrel{\sigma}{\mapsto}(\sigma S,\Pi\sigma^{-1})$, for
$\sigma\in\cs_n$.  The {\bf relative $n$-resonances} are defined to be
the orbits of this $\cs_n$-action. We let $[S,\Pi]$ denote the
relative $n$-resonance represented by the relative $n$-cut~$(S,\Pi)$.
\end{df}

When $S\in\catr_n$ and $\Pi\subseteq P(n)$, $\Pi$ is $S$-closed, it is
convenient to use the notation $Q(S,\Pi)$ to denote the relative cut
$(S\setminus\Pi,\Pi)$. Clearly we have $(S,\Pi)=Q(\tspan S,\Pi)$.
Analogously, $[Q(S,\Pi)]$ denotes the relative resonance
$[S\setminus\Pi,\Pi]$. We use these two notations interchangeably
depending on which one is more natural in the current context.

The special case of the particular importance for our computations in
the later sections is that of $Q(S,\pi\dar S)$, where $\pi$ is
a~partition of $[n]$ with $m$ parts. In this case, we call
$(S\setminus(\pi\dar S),\pi\dar S)$ the relative $(n,m)$-cut
associated to $S$ and $\pi$.

By the Definition~\ref{dfrel}, the relative cut $(S,\Pi)=((\lspan
S)\setminus\Pi,\Pi)$ consists of two parts. We intuitively think of
$(\lspan S)\setminus\Pi$ as the set of all resonances which survive
the shrinking of the strata associated to the elements of $\Pi$, so it
is natural to call them {\it surviving elements}. We also think of
$\Pi$ as the set of all partitions whose associated strata are shrunk
to the infinity point, so, accordingly, we call them {\it partitions
  at infinity}.

\subsection{Direct products of relative resonances.} $\,$

\begin{df} \label{df3.3}$\,$
For relative resonances $(S,\Pi)$ and $(T,\Lambda)$ we define
$$(S,\Pi)\times (T,\Lambda)=
(S\times T,(\Pi\times P(m))\cup (P(n)\times\Lambda)).$$
Clearly the orbit $[(S,\Pi)\times (T,\Lambda)]$ does not depend
on the choice of representatives of the orbits $[S,\Pi]$ and
$[T,\Lambda]$, so we may define $[S,\Pi]\times[T,\Lambda]$ to be
$[(S,\Pi)\times (T,\Lambda)]$.
\end{df}

\nin
The following special cases are of particular importance for our
computation:

\nin
(1) {\bf A direct product of two resonances.}

\nin For an~$m$-cut~$S$, and an~$n$-cut~$T$, we have 
$$S\times T=\{(x_1,\dots,x_m,y_1,\dots,y_n)\,|\,(x_1,\dots,x_m)\in S,
(y_1,\dots,y_n)\in T\}\in\catr_{m+n},$$ 
and $[S]\times[T]=[S\times T]$.

\nin
(2) {\bf A direct product of a relative resonance and a resonance.}

\nin For $S\in\catr_n$, $\Pi\subseteq P(n)$ an $S$-closed set of
partitions, and $T\in\catr_k$, we have $Q(S,\Pi)\times T=Q(S\times
T,\widetilde\Pi)$, where $\widetilde\Pi=\Pi\times
P(\{n+1,n+2,\dots,n+k\})$, and $[Q(S,\Pi)]\times [T]=[Q(S,\Pi)\times T]$.

\vskip4pt   

\nin {\bf Example.}
\begin{multline*}
Q(\{(0,0,0),\pm(1,-1,-1),\pm(0,1,-1)\},\{1\}\{23\})=\\
\{(0)\}\times Q(\{(0,0),\pm(1,-1)\},\{12\}).
\end{multline*}

\begin{rem}
One can define a~category, called {\bf relative resonance category}, 
whose set of objects is the set of all relative $n$-cuts. 
A~new structure which it has in comparison to~$\catr$ is provided by 
``shrinking morphisms'': $(S,\Pi)\rsa(T,\Lambda)$, for 
$S,T\subseteq\{-1,0,1\}^n$, $P(n)\supseteq\Lambda\supseteq\Pi$, such 
that $(\lspan S)\setminus\Lambda=T$. They correspond to shrinking 
strata to infinity. 
\end{rem}

\subsection{Resonance Functors} $\,$

Given a~functor $\cf:\catr\lra\ctop^*$, we introduce the following
notation:
$$\cf(Q(S,\Pi))=\cf(S)\bigg/\bigcup_{\un(\pi)\in\Pi}\im\cf(\pi S
\stackrel{\pi}{\hra}S).$$

\begin{df} \label{dfresf}
  A functor $\cf:\catr\lra\ctop^*$ is called a~{\bf resonance functor}
if it satisfies the following axioms:
\begin{enumerate}
\item [(A1)] {\bf Inclusion axiom.} 
  
  \nin If $S\in\catr_n$, and $\pi\in OP(n)$, then $\cf(\pi S
\stackrel{\pi}{\hra}S)$ is an~inclusion map, and 
$\cf(S)\big/\im\cf(\pi S\stackrel{\pi}{\hra}S)\simeq\cf(Q(S,\pi\dar S))$.

\item [(A2)] {\bf Relative resonance axiom.}

  \nin If, for some $S,T\in\catr_n$, and $\Pi,\Lambda\subseteq P(n)$,
  $[Q(S,\Pi)]=[Q(T,\Lambda)]$, then
  $\cf(Q(S,\Pi))\simeq\cf(Q(T,\Lambda))$.

\item [(A3)] {\bf Direct product axiom.}

\nin For two relative $n$-cuts $(S,\Pi)$ and $(T,\Lambda)$ we have
  $$\cf(S,\Pi)\times\cf(T,\Lambda)\simeq\cf((S,\Pi)\times(T,\Lambda)).$$

\end{enumerate}
\end{df}

Given $S\in\catr_n$, and $\pi\in OP(n)$, let $i_{S,\pi}$ denote the
inclusion map $\cf(\pi S\stackrel{\pi}{\hra}S)$. There is a~canonical
homology long exact sequence associated to the triple
\begin{equation} \label{triple}
 \cf(\pi S)\stackrel{i_{S,\pi}}{\hra}\cf(S)\stackrel{p}{\lra}
\cf(Q(S,\pi\dar S)),
\end{equation}
namely
\begin{multline} \label{stls}
\dots\stackrel{\partial_*}{\lra}\widetilde H_n(\cf(\pi S))
\stackrel{(i_{S,\pi})_*}{\lra}\widetilde H_n(\cf(S))\stackrel{p_*}{\lra}
\widetilde H_n(\cf(Q(S,\pi\dar S)))\stackrel{\partial_*}{\lra}  \\
\widetilde H_{n-1}(\cf(\pi S))\stackrel{(i_{S,\pi})_*}{\lra}\dots
\end{multline}
We call \eqref{triple}, resp.~\eqref{stls}, the {\it standard triple},
resp.~the {\it standard long exact sequence} associated to the
morphism $\pi S\stackrel{\pi}{\hra}S$ and the functor $\cf$ (usually
$\cf$ is fixed, so its mentioning is omitted).

     \section{First applications}

\subsection{Resonance compatible compactifications.} $\,$
\label{ss4.1}

As mentioned in the introduction we shall now look at the natural
strata of the spaces $X^{(n)}$. The strata are defined by point
coincidences and are indexed by number partitions of $n$. Let
$\Sigma_\lambda^X$ denote the stratum indexed by $\lambda$.

Let $\lambda$ be a number partition of $n$ and let $\tilde\lambda$ be
$\lambda$ with some fixed order on the parts. Then $\tilde\lambda$ can
be thought of as a~vector with positive integer coordinates
in~$\br^n$. Let $S_{\tilde\lambda}$ be the set
$\{x\in\{-1,0,1\}^n\,|\,\langle x,\tilde\lambda\rangle =0\}$.
Obviously, $S_{\tilde\lambda}$ is an~$n$-cut and the $n$-resonance
$S_\lambda$, which it defines, does not depend on the choice of
$\ti\lambda$, but only on the number partition~$\lambda$.

The crucial topological observation is that if $\nu$ is another
partition of $n$, such that $S_\lambda=S_\nu$, then the spaces $\slx$
and $\Sigma_\nu^X$ are homeomorphic. This is precisely the fact which
leads one to introduce resonances and the surrounding combinatorial
framework and to forget about the number partitions themselves.

This allows us to introduce a~functor $\cf$ mapping
$S_{\tilde\lambda}$ to $\Sigma_\lambda^X$; the morphisms map
accordingly. Clearly, $\cf(1^l)=X^{(l)}$. One can observe in this
example the justification for the names which we chose for the
morphisms of $\catr$: ``inclusions'' and ``gluings''. Furthermore, it
is easy, in this case, to verify the axioms of the
Definition~\ref{dfresf}, and hence to conclude that $\cf$ is
a~resonance functor. The only nontrivial point is the verification of
the second part of (A1), which we do in the next proposition.

\begin{prop} \label{prop4.1}
  Let $S$ be an $n$-cut and $\pi\in OP(n)$. Then $\cf(\nu
  S)\subseteq\cf(\pi S)$ if and only if $\un(\nu)\in\un(\pi)\dar S$.
\end{prop}
\pr It is obvious that all the steps of the definition of
$\un(\pi)\dar S$ which change the partition preserve the property
$\cf(\nu S)\subseteq\cf(\pi S)$, hence the {\it if} direction follows.

Assume now $\cf(\nu S)\subseteq\cf(\pi S)$. This means that there
exists $\tau\in OP(m)$, where $m$ is the number of parts of $\pi$,
such that $\cf(\tau\pi S)=\cf(\nu S)$. By definition,
$\tau\circ\pi\in\un(\pi)\dar S$. Now, we can reach $\un(\nu)$ from
$\un(\tau\circ\pi)$ by moves of type (3) from the definition of the
relative resonances. 

Indeed, if $\cf(\tau\pi S)=\cf(\nu S)=\Sigma_\lambda^X$, then the
sizes of the resulting blocks after gluing along $\tau\circ\pi$ and
along $\nu$ are the same. For every block $b$ of $\lambda$ we can go,
by means of moves of type (3), from the block of $\un(\tau\circ\pi)$
which glues to $b$ to the block of $\un(\nu)$ which glues to $b$.
Since we can do it for any block of $\lambda$, we can go from
$\un(\tau\circ\pi)$ to $\un(\nu)$, and hence 
$\un(\nu)\in\un(\pi)\dar S$. 
\qed

\vskip3pt

\nin In the context of this stratification the following central
question arises.

\vskip3pt

\nin {\bf The Main Problem.} {\it (Arnold, Shapiro,~\cite{Sh}).
  Describe an algorithm which, for a~given resonance $\lambda$, would
  compute the Betti numbers of $\Sigma_\lambda^{S^1}$, or
  $\Sigma_\lambda^{S^2}$.}

\vskip3pt

The case of the strata $\Sigma_\lambda^{S^1}$ is simpler, essentially
because of the following elementary, but important property of smash
products: if $X$ and $Y$ are pointed spaces and $X$ is contractible,
then $X\wedge Y$ is also contractible.

In the subsequent subsections we shall look at a~few interesting
special cases, and also will be able to say a~few things about the
general problem.

\subsection{Resonances $(a^k,1^l)$.} $\,$ \label{ss4.2}

Let $a,k,l$ be positive integers such that $a\geq 2$. Let $S$ be the
$(l+k)$-cut consisting of all the elements of
$\{-1,0,1\}^{l+k}$, which are orthogonal to the vector
$(\underbrace{1,\dots,1}_l,\underbrace{a,\dots,a}_k)$. Clearly,
the~$(l+k)$-resonance $[S]$ is equal to $(a^k,1^l)$. The case $l<a$ is
not very interesting, since then $(a^k,1^l)=(1^k)\times(1^l)$.
Therefore we may assume that $l\geq a$.

We would like to understand the topological properties of the space
$\cf(a^k,1^l)$. In general, this is rather hard. However, as the
following theorem shows, it is possible under some additional
conditions on $\cf$.

\begin{thm} \label{thm4.2}
  Let $\cf:\catr\lra\ctop^*$ be a~resonance functor, such that
  $\cf(1^l)$ is contractible for $l\geq 2$. Let $l=am+\epsilon$, where
  $0\leq\epsilon\leq a-1$.
  \begin{enumerate}
  \item [(a)] If $k\neq 1$, or $\epsilon\geq 2$, then $\cf(a^k,1^l)$
    is contractible.
  \item [(b)] If $k=1$, and $\epsilon\in\{0,1\}$, then
  \begin{equation} \label{eq:4.1}
   \cf(a^k,1^l)\simeq\susp^m(\cf(1)^{m+\epsilon+1}), 
  \end{equation}
  where $\cf(1)^{m+\epsilon+1}$ denotes the $(m+\epsilon+1)$-fold
  smash product.
  \end{enumerate}
\end{thm}

Since for the resonance functor $\cf$ described in the
subsection~\ref{ss4.1} we have $\cf(1^l)=X^{(l)}$, we have the
following corollary.

\begin{crl} \label{crl4.3}
  If $X^{(l)}$ is contractible for $l\geq 2$, then
  \begin{enumerate}
  \item [(a)] If $k\neq 1$, or $\epsilon\geq 2$ (again
    $l=am+\epsilon$), then $\Sigma_{(a^k,1^l)}^X$ is contractible.
  \item [(b)] If $k=1$, and $\epsilon\in\{0,1\}$, then
    $\Sigma_{(a^k,1^l)}^X\simeq\susp^m(X^{m+\epsilon+1})$, where
    $X^{m+\epsilon+1}$ denotes the $(m+\epsilon+1)$-fold smash
    product.
  \end{enumerate}
\end{crl}

\nin {\bf Note.} Clearly, $(S^1)^{(l)}$ is contractible for $l\geq 2$,
so the Corollary~\ref{crl4.3} is valid. In this situation, the case
$k>1$ was proved in~\cite{Ko2}, and the case $k=1$ in~\cite{BW,SW}.

\vskip4pt
Before we proceed with proving the Theorem~\ref{thm4.2} we need
a~crucial lemma. Let $\pi\in P(k+l)$ be 
$(\{1,\dots,a\},\{a+1\},\{a+2\},\dots,\{k+l\})$. It is immediate
that $[\tilde\pi S]=(a^{k+1},1^{l-a})$, if $\un(\tilde\pi)=\pi$.

\begin{lm} \label{lm4.4}
  Let $S$ be as above, $T\in\catr_l$, such that $[T]=(1^l)$, and let
  $\nu$ be the partition
  $(\{1,\dots,a\},\{a+1\},\{a+2\},\dots,\{l\})$, then we have
  \begin{equation}
    \label{eq:4.2}
    [Q(S,\pi\dar S)]=[Q(T,\nu\dar T)]\times(1^k).
  \end{equation}
\end{lm}

\nin {\bf Note.} Lemma~\ref{lm4.4} is a~special case of the
Lemma~\ref{lm4.6}, however we choose to include a~separate proof for
it for two reasons: firstly, it is the first, still not too technical
example of investigating the combinatorial structure of the resonance
category, which is a~new object; secondly, the particular case of
$(a^k,1^l)$ resonances was a~subject of substantial previous
attention.

\vskip3pt
\nin {\bf Proof of the Lemma~\ref{lm4.4}.} 

\nin Recall that by the definition of the direct product,
$$[Q(T,\nu\dar T)]\times(1^k)=[Q(T\times U,(\nu\dar T)\times
P(\{l+1,\dots,l+k\}))],$$ 
where $U\in\catr_k$ and $[U]=(1^k)$. Clearly, $(\nu\dar T)\times
P(\{l+1,\dots,l+k\})=\pi\dar S$, hence we just need to show that
$S\setminus(\pi\dar S)=(T\times U)\setminus((\nu\dar T)\times 
P(\{l+1,\dots,l+k\}))$. Note that $(T\times U)\setminus((\nu\dar T)\times 
P(\{l+1,\dots,l+k\}))=(T\setminus(\nu\dar T))\times U$.
Furthermore, 
$$S=\bigg\{(x_1,\dots,x_{l+k})\in\{-1,0,1\}^{l+k}\,\Bigm|\,
\sum_{j=l+1}^{l+k}x_j+a\sum_{i=1}^{l}x_i=0\bigg\},$$
and the set which we need to remove from $S$ to get 
$S\setminus(\pi\dar S)$ is
\begin{multline*}
\bigg\{(x_1,\dots,x_{l+k})\in\{-1,0,1\}^{l+k}
\,\Bigm|\,\sum_{j=l+1}^{l+k}x_j+a\sum_{i=1}^{l}x_i=0, \\
\max(|\plus(x_1,\dots,x_l)|,|\minus(x_1,\dots,x_l)|)\geq a\bigg\}.
\end{multline*}
Therefore, by the definition of the relative resonances, we have
\begin{multline*}
S\setminus(\pi\dar S)=\bigg\{(x_1,\dots,x_{l+k})\in\{-1,0,1\}^{l+k}\,\Bigm|\, \\
\sum_{i=1}^l x_i=0,\,\,\sum_{j=l+1}^{l+k}x_j=0,\,\,
|\plus(x_1,\dots,x_l)|<a\bigg\}.
\end{multline*}
On the other hand, 
$(1^k)=[\{(y_1,\dots,y_k)\in\{-1,0,1\}^k\,|\,\sum_{j=1}^k y_j=0\}]$,
and
$$T\setminus(\nu\dar T)=\bigg\{(z_1,\dots,z_l)\in\{-1,0,1\}^l\,\Bigm|\,
\sum_{i=1}^l z_i=0,\,\,|\plus(z_1,\dots,z_l)|<a\bigg\},$$
which proves~\eqref{eq:4.2}. \qed

\vskip3pt
\nin {\bf Proof of the Theorem~\ref{thm4.2}.} 

\nin {\bf (a)} We use induction on~$l$. The case $l<a$ can be taken as
an~induction base, since then $(a^k,1^l)=(1^k)\times(1^l)$, hence, by
the axiom (A3), $\cf(a^k,1^l)=\cf(1^k)\wedge\cf(1^l)$, which is
contractible, since $\cf(1^k)$ is. Thus we assume that $l\geq a$, and
$\cf(a^k,1^{l'})$ is contractible for all $l'<l$.

Let $S$ and $\pi$ be as in the Lemma~\ref{lm4.4}. The standard triple
associated to the morphism $\pi S\stackrel{\pi}{\hra}S$ is
$\cf(a^{k+1},1^{l-a})\hra\cf(a^k,1^l)\ra\cf(a^k,1^l)/\cf(a^{k+1},1^{l-a})$.
Since, by the induction assumption, $\cf(a^{k+1},1^{l-a})$ is
contractible, we conclude that $\cf(a^k,1^l)\simeq\cf(a^k,1^l)/
\cf(a^{k+1},1^{l-a})$.

Basically by the definition, we have
$\cf(a^k,1^l)/\cf(a^{k+1},1^{l-a})= \cf(Q(S,\pi\dar S))$. On the other
hand, we have proved in the Lemma~\ref{lm4.4} that $[Q(S,\pi\dar
S)]=Q(T,\nu\dar T)\times (1^k)$, where $T$ and $\nu$ are described in
the formulation of that lemma. By axioms (A2) and (A3) we get that
$\cf(Q(S,\pi\dar S))\simeq\cf(Q(T,\nu\dar T))\wedge \cf(1^k)$, which
is contractible, since $\cf(1^k)$ is. Therefore, $\cf(a^k,1^l)$ is
also contractible.

\nin {\bf (b)} The argument is very similar to (a). We again assume
$l\geq a$, which implies $l\geq 2$. By the using the same ordered set
partition~$\pi$ as in (a), we get that $\cf(a,1^l)\simeq
\cf(a,1^l)/\cf(a^2,1^{l-a})$. Further, by Lemma~\ref{lm4.4} and the
axioms (A2) and (A3) we conclude that $\cf(a,1^l)\simeq
\cf(1)\wedge(\cf(1^l)/\cf(a,1^{l-a}))$. Since $\cf(1^l)$ is
contractible, we get 
\begin{equation}
  \label{eq:4.3}
  \cf(a,1^l)\simeq \cf(1)\wedge \susp\cf(a,1^{l-a}).
\end{equation}
Since $\cf(a)=\cf(1)$, $\cf(a,1)=\cf(1)\wedge\cf(1)$, and $\cf(a,1^l)$
is contractible if $2\leq l<a$, we obtain \eqref{eq:4.1} by the
repeated usage of~\eqref{eq:4.3}. \qed

\subsection{Resonances $(a^k,b^l)$.} $\,$

The algebraic invariants of these strata have not been computed
before, not even in the case $X=S^1$, and $\cf$ - the standard
resonance functor associated to the stratification of $X^{(n)}$.

We would like to apply a~technique similar to the one used in the
subsection~\ref{ss4.2}. A~problem is that, once one starts to ``glue''
$a$'s, one cannot get $b$'s in the same way as one could in the
previous section from 1's. Thus, we are forced to consider a~more 
general case of resonances, namely $(g^m,a^k,b^l)$, where $g$ is
the least common multiple of $a$ and $b$. Assume $g=a\cdot\bar a=
b\cdot\bar b$, and $b>a\geq 2$. Analogously with the Theorem~\ref{thm4.2}
we have the following result.

\begin{thm} \label{thm4.5}
  Let $\cf$ be as in the Theorem~\ref{thm4.2}. Let furthermore 
$k=x\cdot\bar a+\epsilon_1$, $l=y\cdot\bar b+\epsilon_2$, where 
$0\leq\epsilon_1<\bar a$, $0\leq\epsilon_2<\bar b$. Then
\begin{equation}
  \label{eq:4.4}
  \cf(g^m,a^k,b^l)\simeq
     \begin{cases}
        \susp^{x+y+m-1}(\cf(1)^{x+y+m+\epsilon_1+\epsilon_2}), & 
        \text{ if } m,\epsilon_1,\epsilon_2\in\{0,1\}; \\
        \text{ point, } & \text{ otherwise. } 
     \end{cases} 
\end{equation}
\end{thm}

Just as in the subsection~\ref{ss4.2} (Corollary~\ref{crl4.3}), the
Theorem~\ref{thm4.5} is true if one replaces $\cf(\lambda)$ with
$\Sigma_\lambda^{S^1}$.

The proof of the Theorem~\ref{thm4.5} follows the same general scheme
as that of the Theorem~\ref{thm4.2}, but the technical details are
more numerous. Again there is a~crucial combinatorial lemma.

Let $S$ be an~$(m+k+l)$-cut consisting of all the elements of
$\{-1,0,1\}^{m+k+l}$ which are orthogonal to the vector
$(\underbrace{a,\dots,a}_k,\underbrace{b,\dots,b}_l,
\underbrace{g,\dots,g}_m)$. Assume $k\geq\bar a$, and let an~unordered
set partition~$\pi$ be equal to $(\{1,\dots,\bar a\},\{\bar
a+1\},\dots,\{k+l+m\})$. We see that $[S]=(g^m,a^k,b^l)$, and
$[\tilde\pi S]=(g^{m+1},a^{k-\bar a},b^l)$, if $\pi=\un(\tilde\pi)$.

\begin{lm} \label{lm4.6}
Let $T\in\catr_k$, such that $[T]=(1^k)$, and 
$\nu=(\{1,\dots,\bar a\},\{\bar a+1\},\dots,\{k\})$, then
\begin{equation}
  \label{eq:4.5}
  [Q(S,\pi\dar S)]=[Q(T,\nu\dar T)]\times(\bar b^m,1^l).
\end{equation}
\end{lm}
\pr Again, it is easy to see that the sets of the partitions at infinity
on both sides of~\eqref{eq:4.5} coincide. Indeed,
$$[Q(T,\nu)]\times(\bar b^m,1^l)=[Q(T\times U,(\nu\dar T)\times
P(\{k+1,\dots,k+m+l\}))],$$
where $U\in\catr_{m+l}$, such that
$[U]=({\bar b}^m,1^l)$, and $(\nu\dar T)\times
P(\{k+1,\dots,k+m+l\})=\pi\dar S$. Also, we again have the equality
$$(T\times U)\setminus((\nu\dar T)\times
P(\{k+1,\dots,k+m+l\}))=T\setminus(\nu\dar T)\times U,$$
which greatly
helps to prove that the sets if the surviving elements on the two
sides of~\eqref{eq:4.5} coincide. 

By the definition
$$S=\bigg\{(x_1,\dots,x_{k+l+m})\in\{-1,0,1\}^{k+l+m}\,\Bigm|\,
a\sum_{i=1}^{k}x_i+b\sum_{i=k+1}^{k+l}x_i+g\sum_{i=k+l+1}^{k+l+m}x_i=0\bigg\},$$
and, again, the set which we have to remove from $S$ to get 
$S\setminus(\pi\dar S)$ is
\begin{multline*}
\bigg\{(x_1,\dots,x_{k+l+m})\in\{-1,0,1\}^{k+l+m}\,\Bigm|\,
a\sum_{i=1}^{k}x_i+b\sum_{i=k+1}^{k+l}x_i+ \\ g\sum_{i=k+l+1}^{k+l+m}x_i=0,\,\,
\max(|\plus(x_1,\dots,x_k)|,|\minus(x_1,\dots,x_k)|)\geq\bar a\bigg\}.$$
\end{multline*}
By the definition of the relative resonances and some elementary number 
theory we conclude that 
\begin{multline*}
S\setminus(\pi\dar S)=\bigg\{(x_1,\dots,x_{k+l+m})\in\{-1,0,1\}^{k+l+m}\,\Bigm|\,
 |\plus(x_1,\dots,x_k)|<\bar a,  \\
\sum_{i=1}^{k}x_i=0,\quad
b\sum_{i=k+1}^{k+l}x_i+g\sum_{i=k+l+1}^{k+l+m}x_i=0\bigg\}.$$
\end{multline*}
The number theory argument which we need is that if $ax+by+\lcm(a,b)z=0$,
then $\bar a\,\big|\,x$, where $\bar a\cdot a=\lcm(a,b)$. This can be 
seen by, for example, noticing that if $ax+by+\lcm(a,b)z=0$, then  
$b\,\big|\,ax$, but since also $a\,\big|\,ax$, we have 
$\lcm(a,b)\,\big|\,ax$, hence $\bar a\,\big|\,x$.

The equation~\eqref{eq:4.5} follows now from the earlier observations
together with the equalities
$$T\setminus(\nu\dar T)=\bigg\{(x_1,\dots,x_k)\in\{-1,0,1\}^k\,\Bigm|\,
|\plus(x_1,\dots,x_k)|<\bar a,\,\,\sum_{i=1}^{k}x_i=0\bigg\},$$
and 
$$(\bar b^m,1^l)=\bigg[\bigg\{(y_1,\dots,y_{m+l})\in\{-1,0,1\}^{m+l}
\,\Bigm|\,\sum_{i=1}^l y_i+\bar b\sum_{i=l+1}^{l+m}x_i=0\bigg\}\bigg].\qed$$

\vskip3pt
\nin {\bf Proof of the Theorem~\ref{thm4.5}.}
  The cases $k<\bar a$ and $l<\bar b$ are easily reduced to the 
Theorem~\ref{thm4.2}. Assume therefore that $k\geq\bar a$ and 
$l\geq\bar b$. Recall also that $b>a\geq 2$, and hence $\bar a\geq 2$. 

Let $S$ and $\pi$ be as in the formulation of the Lemma~\ref{lm4.6}. 
The standard triple associated to the morphism 
$\pi S\stackrel{\pi}{\hra}S$ is
\begin{equation}
  \label{eq:4.6}
  \cf(g^{m+1},a^{k-\bar a},b^l)\hra\cf(g^m,a^k,b^l)\ra
  \cf(g^m,a^k,b^l)/\cf(g^{m+1},a^{k-\bar a},b^l).
\end{equation}
We break the rest of the proof into 3 cases.

\nin {\bf Case $m\geq 2$.} Again, we prove that $\cf(g^m,a^k,b^l)$ is
contractible by induction on~$k$. This is clear if $k<\bar a$. If
$k\geq\bar a$, it follows from~\eqref{eq:4.6} that $\cf(g^m,a^k,b^l)
\simeq\cf(g^m,a^k,b^l)/\cf(g^{m+1},a^{k-\bar a},b^l)=\cf(Q(S,\pi\dar
S))$. By Lemma~\ref{lm4.6} we conclude that $\cf(g^m,a^k,b^l)\simeq
\cf(Q(T,\nu\dar T))\wedge\cf(\bar b^m,1^l)$. By the
Theorem~\ref{thm4.2}, $\cf(\bar b^m,1^l)$ is contractible, hence so is
$\cf(g^m,a^k,b^l)$.

\nin {\bf Case $m=0$.} By Lemma~\ref{lm4.6} we get that
$\cf(Q(S,\pi\dar S))\simeq\cf(Q(T,\nu\dar T))\wedge\cf(1^l)$. Since
$l\geq 2$, we have that $\cf(1^l)$ is contractible, hence so is
$\cf(Q(S,\pi\dar S))= \cf(a^k,b^l)/\cf(g,a^{k-\bar a},b^l)$.
Therefore, by~\eqref{eq:4.6} $\cf(a^k,b^l)\simeq\cf(g,a^{k-\bar
  a},b^l)$.

\nin {\bf Case $m=1$.} Since $\cf(g^2,a^{k-\bar a},b^l)$ is
contractible, we conclude by~\eqref{eq:4.6} that
$\cf(g,a^k,b^l)\simeq\cf(g,a^k,b^l)/ \cf(g^2,a^{k-\bar
  a},b^l)=\cf(Q(S,\pi\dar S))$. By Lemma~\ref{lm4.6}, and the
properties of the resonance functors, we have
\begin{equation}
  \label{eq:4.7}
\cf(g,a^k,b^l)\simeq\cf(\bar b,1^l)\wedge(\cf(1^k)/\cf(\bar a,1^{k-\bar a}))
\simeq\cf(\bar b,1^l)\wedge\susp(\cf(\bar a,1^{k-\bar a})).
\end{equation} 
By the repeated usage of~\eqref{eq:4.7} we obtain~\eqref{eq:4.4}.
\qed

     \section{Sequential resonances.} 

\subsection{The structure theory of strata associated to sequential 
resonances.} $\,$

\begin{df}\label{dfseqres}
  Let $\lambda=(\lambda_1,\dots,\lambda_n)$,
  $\lambda_1\leq\dots\leq\lambda_n$, be a~number partition.  We call
  $\lambda$ {\bf sequential} if, whenever $\sum_{i\in
    I}\lambda_i=\sum_{j\in J}\lambda_j$, and $q\in I$, such that
  $q=\max(I\cup J)$, then there exists $\widetilde J\subseteq J$, such
  that $\lambda_q=\sum_{j\in\widetilde J}\lambda_j$.
  
  Correspondingly, we call a~resonance $S$ sequential, if it can be
  associated to a~sequential partition.
\end{df}

\nin Note that the set of sequential partitions is closed under
removing blocks.

\vskip3pt

\nin {\bf Examples of sequential partitions:}
\begin{enumerate}
\item all partitions whose blocks are equal to powers of some number;
\item $(a^k,b^l,1^m)$, such that $a>bl$; more generally
  $(a_1^{k_1},\dots,a_t^{k_t},1^m)$, such that $a_i>\sum_{j=i+1}^t a_j
  k_j$, for all $i\in[t]$.
\end{enumerate}

Through the rest of this subsection, we let $\lambda$ be as in the
Definition~\ref{dfseqres}. For such $\lambda$ we use the following
additional notations:
\begin{itemize}
\item $mm(\lambda)=|\{i\in[n]\,|\,\lambda_i=\lambda_n\}|$. 
In other words
$\lambda_{n-mm(\lambda)}\neq\lambda_{n-mm(\lambda)+1}=\dots=\lambda_n$.
\item $I(\lambda)\subseteq[n]$ is the lexicographically maximal set
  (see below the convention that we use to order lexicographically),
  such that $|I(\lambda)|\geq 2$, and $\lambda_n=\sum_{i\in
    I(\lambda)}\lambda_i$. Note that it may happen that $I(\lambda)$
  does not exist, in which case $\cf(\lambda)\simeq
  \cf(\lambda_1,\dots,\lambda_{n-mm(\lambda)})\wedge\cf(1^{mm(\lambda)})$,
  and can be dealt with by induction.
\end{itemize}

Let $n$ be a~positive integer. We use the following convention for the
lexicographic order on $[n]$. For $A=\{a_1,\dots,a_k\}$,
$B=\{b_1,\dots,b_m\}$, $A,B\subseteq[n]$, $a_1\leq\dots\leq a_k$,
$b_1\leq\dots\leq b_m$, we say that $A$ is lexicographically larger
than~$B$ if, either $A\supseteq B$, or there exists $q<\min(k,m)$,
such that $a_k=b_m$, $a_{k-1}=b_{m-1}$, $\dots$,
$a_{k-q+1}=b_{m-q+1}$, and $a_{k-q}>b_{m-q}$.

\begin{prop} \label{prop4.9}
  If $\lambda=(\lambda_1,\dots,\lambda_n)$,
  $\lambda_1\leq\dots\leq\lambda_n$, is a sequential partition, then
  so is $\bar\lambda=(\lambda_{j_1},\dots,\lambda_{j_t},\sum_{i\in
    I(\lambda)}\lambda_i)$, where $t=n-|I(\lambda)|$, and
  $\{j_1,\dots,j_t\}=[n]\setminus I(\lambda)$.
\end{prop}
\pr Let $\bar\lambda_1=\lambda_{j_1},\dots,\bar\lambda_t=\lambda_{j_t}$,
$\bar\lambda_{t+1}=\sum_{i\in I(\lambda)}\lambda_i$.
We need to check the condition of the Definition~\ref{dfseqres}
for the identity
\begin{equation}
  \label{eq:4.8a}
  \sum_{i\in I}\bar\lambda_i=\sum_{j\in J}\bar\lambda_j.
\end{equation}

If $t+1\not\in I\cup J$, then it follows from the assumption that
$\lambda$ is sequential. Assume $t+1\in I$. If
$\bar\lambda_j=\lambda_n$, for some $j\in J$, take $\widetilde
J=\{j\}$, and we are done. If $\bar\lambda_i=\lambda_n$, for some
$i\in I\setminus\{t+1\}$, then, since $\lambda$ is sequential, there
exists $\widetilde J\subseteq J$, such that $\sum_{j\in\widetilde
  J}\bar\lambda_j=\lambda_n=\bar\lambda_{t+1}$, and we are done again.

Finally, assume $\bar\lambda_i\neq\lambda_n$, for $i\in(I\cup
J)\setminus\{t+1\}$. Substituting $\lambda_n$ instead of
$\bar\lambda_{t+1}$ into the identity~\eqref{eq:4.8a} is allowed,
since $\lambda_n$ does not appear among $\{\bar\lambda_i\}_{i\in(I\cup
  J)\setminus\{t+1\}}$. This gives us an identity for $\lambda$, and
again, since $\lambda$ is sequential, we find the desired set
$\widetilde J\subseteq J$, such that $\sum_{j\in\widetilde
  J}\bar\lambda_j=\bar\lambda_{t+1}$. \qed

\vskip4pt

Let $S\in\catr_n$ be the set of all elements of $\{-1,0,1\}^n$, which
are orthogonal to the vector $\lambda=(\lambda_1,\dots,\lambda_n)$.
Clearly, $[S]=\lambda$. Let $\pi\in P(n)$ be the partition whose only
nonsingleton block is given by $I(\lambda)$. The next lemma expresses
the main combinatorial property of sequential partitions.

\begin{lm}\label{green}
  Let $\tau\in P(n)$ be a~partition which has only one nonsingleton
block~$B$, and assume $\lambda_n=\sum_{i\in B}\lambda_i$. 
Then $\tau\in\pi\dar S$.
\end{lm}
\pr Assume there exists partitions $\tau$ as in the formulation of the
lemma, such that $\tau\not\in\pi\dar S$. Choose one so that the
block~$B$ is lexicographically largest possible. Let $C=B\cap
I(\lambda)$. By the definition of $I(\lambda)$, and the choice of $B$,
we have $\sum_{i\in I(\lambda)\setminus C}\lambda_i=
\sum_{j\in B\setminus C}\lambda_j$, and $q\in I(\lambda)\setminus C$,
where $q=\max((I(\lambda)\cup B)\setminus C)$.

Since partition $\lambda$ is sequential, there exists $D\subseteq
B\setminus C$, such that $\lambda_q=\sum_{j\in D}\lambda_j$.  Let
$\gamma\in P(n)$ be the partition whose only nonsingleton block is
$G=(B\setminus D)\cup\{q\}$. Clearly, $\sum_{i\in G}\lambda_i=
\lambda_n$, and $|G|\geq 2$. By the choice of $q$, $G$ is 
lexicographically larger than $B$, hence $\gamma\in\pi\dar S$.

Let furthermore $\tilde\gamma\in P(n)$ be the partition having two
nonsingleton blocks: $D$ and $G$. By the Definition~\ref{df3.1}(2) if
$\gamma\in\pi\dar S$, then $\tilde\gamma\in\pi\dar S$. By the
Definition~\ref{df3.1}(3), if $\tilde\gamma\in\pi\dar S$, then
$\tau\in\pi\dar S$, which yields a~contradiction.
\qed

\vskip4pt

Let $T\in\catr_{n-mm(\lambda)}$ be the set of all elements of
$\{-1,0,1\}^{n-mm(\lambda)}$, which are orthogonal to the vector
$(\lambda_1,\dots,\lambda_{n-mm(\lambda)})$. Let $\nu\in
P(n-mm(\lambda))$ be the partition whose only nonsingleton block is
given by $I(\lambda)$. We are now ready to state the combinatorial
result which is crucial for our topological applications.

\begin{lm} \label{lm4.11}
  \begin{equation}
    \label{eq:4.9}
    [Q(S,\pi\dar S)]=[Q(T,\nu\dar T)]\times(1^{mm(\lambda)}).
  \end{equation}
\end{lm}
\pr By definition we must verify that the sets of partitions at
infinity and the surviving elements coincide on both sides of the
equation~\eqref{eq:4.9}. 

Let us start with the partitions at infinity. Filtered through the
Proposition~\ref{prop4.1}, the identity $\pi\dar S=(\nu\dar T)\times
P(\{n-mm(\lambda)+1,\dots,n\})$ becomes essentially tautological. Both
sides consist of the partitions $\tau=(\tau_1,\dots,\tau_k)\in P(n)$,
such that the number partition 
$(\sum_{i\in\tau_1}\lambda_i,\dots,\sum_{i\in\tau_k}\lambda_i)$ 
can be obtained from the number partition 
$(\lambda_{j_1},\dots,\lambda_{j_t},\sum_{i\in I(\lambda)}\lambda_i)$,
where $\{j_1,\dots,j_t\}=[n]\setminus I(\lambda)$, by summing parts.

Let us now look at the surviving elements. Obviously,
$S\setminus(\pi\dar S)\supseteq(T\setminus(\nu\dar T))\times U$, where
$U\in\catr_k$, such that $[U]=(1^{mm(\lambda)})$, and we need to show
the converse inclusion. Let $x=(x_1,\dots,x_n)\in S$, such that
$\sum_{i=n-mm(\lambda)+1}^n x_i\neq 0$ (otherwise
$x\in(T\setminus(\nu\dar T))\times U$), we can assume
$\sum_{i=n-mm(\lambda)+1}^n x_i>0$. Then, since $S$ is a~sequential
resonance, there exists $y=(y_1,\dots,y_n)\in S$, such that
\begin{itemize}
\item if $y_i\neq 0$, then $x_i=y_i$; 
\item $|\plus(y)|=1$, and $\plus(y)\subseteq\{n-mm(\lambda)+1,\dots,n\}$.
\end{itemize}
This, by Lemma~\ref{green}, means that 
$y\not\in S\setminus(\pi\dar S)$, which in turn necessitates 
$x\not\in S\setminus(\pi\dar S)$. 
This finishes the proof of the lemma. 
\qed

\vskip4pt

Just as before, this combinatorial fact about the resonances
translates into a~topological statement, which can be further
strengthened by requiring some additional properties from~$\lambda$.

\begin{df}
  Let $\lambda=(\lambda_1,\dots,\lambda_n)$,
  $\lambda_1\leq\dots\leq\lambda_n$, be a~sequential partition, and
  let $q=\max I(\lambda)$. $\lambda$ is called {\bf strongly
  sequential}, if, either $I(\lambda)$ does not exist, or there exists
  $J\subseteq I(\lambda)\setminus\{q\}$, such that
  $\lambda_q=\sum_{i\in J}\lambda_i$ (note that we do not require
  $|J|\geq 2$).
\end{df}

We are now in a~position to prove the main topological structure
theorem concerning the sequential resonances.

\begin{thm} \label{thm4.13}
  Let $\cf$ be as in the Theorem~\ref{thm4.2}. Let $\lambda$ be 
a~sequential partition, such that $I(\lambda)$ exists, then
\begin{enumerate}
\item if $mm(\lambda)\geq 2$, then $\cf(\lambda)$ is contractible;
\item if $mm(\lambda)=1$, then $\cf(\lambda)\simeq\cf(Q(T,\nu\dar
T))\wedge\cf(1)$, and we have the inclusion triple
$\cf(\mu)\stackrel{i}{\hra}\cf(\lambda_1,\dots,
\lambda_{n-1})\ra\cf(Q(T,\nu\dar T))$, where
$\mu=(\lambda_{j_1},\dots, \lambda_{j_t})$,
$\{j_1,\dots,j_t\}=[n]\setminus I(\lambda)$, and $\nu\in
P(n-mm(\lambda))$ is the partition whose only nonsingleton block is
given by $I(\lambda)$. We set $\cf(\mu)$ to be a~point, if
$I(\lambda)$ does not exist.

If moreover $\lambda$ is strongly sequential, then the map $i$ is
homotopic to a~trivial map (mapping everything to a~point), hence the
triple splits and we conclude that
\begin{equation}
  \label{eq:4.10}
  \cf(\lambda)\simeq(\cf(1)\wedge\cf(\lambda_1,\dots,\lambda_{n-1}))
  \vee\susp(\cf(1)\wedge\cf(\mu)).
\end{equation}
\end{enumerate}
\end{thm}
\pr $\,$

\nin (1) We use induction on $\sum_{i=1}^{n-mm(\lambda)}\lambda_i$.
If $I(\lambda)$ does not exist, then $\lambda_n$ is independent, i.e.,
$\cf(\lambda)\simeq\cf(\lambda_1,\dots,\lambda_{n-mm(\lambda)})
\times\cf(1^{mm(\lambda)})$, and hence $\cf(\lambda)$ is contractible.
Otherwise consider the inclusion triple
\begin{equation}
  \label{eq:4.11}
  \cf(\bar\lambda)\hra\cf(\lambda)\ra\cf(\lambda)/\cf(\bar\lambda)=
\cf(Q(S,\pi\dar S)),
\end{equation}
where $\bar\lambda=(\lambda_{j_1}, \dots,\lambda_{j_t},\sum_{i\in
I(\lambda)}\lambda_i)$, and $\pi\in P(n)$ is the partition whose 
only nonsingleton block is given by $I(\lambda)$. By the induction 
assumption $\cf(\bar\lambda)$ is contractible. On the other hand, 
by Lemma~\ref{lm4.11}, $\cf(Q(S,\pi\dar S))\simeq\cf(Q(T,\nu\dar
T))\wedge\cf(1^{mm(\lambda)})$, which is also contractible if
$mm(\lambda)\geq 2$.

\nin (2) if $mm(\lambda)=1$, then we can conclude from~\eqref{eq:4.11}
that $\cf(\lambda)\simeq\cf(1)\wedge\cf(Q(T,\nu\dar T))$. Next,
consider the inclusion triple 
\begin{equation}
  \label{eq:4.12}
  \cf(\mu)\stackrel{i}{\hra}\cf(\lambda_1,\dots,\lambda_{n-1})\ra
  \cf(Q(T,\nu\dar T)).
\end{equation}
  If $\lambda$ is strongly sequential, then there exists
$J\subseteq I(\lambda)\setminus\{q\}$, such that 
$\lambda_q=\sum_{i\in J}\lambda_i$ (here $q=\max I(\lambda)$). 
The map $i$ factors:
\begin{equation}
  \label{eq:4.13}
  \cf(\mu)\stackrel{i_1}{\hra}\cf(\lambda_{p_1},\dots,
\lambda_{p_{n-1-|J|}},\sum_{i\in I(\lambda)}\lambda_i)
\stackrel{i_2}{\hra}\cf(\lambda_1,\dots,\lambda_{n-1}),
\end{equation}
where $\{p_1,\dots,p_{n-1-|J|}\}=[n-1]\setminus J$. Since
$(\lambda_{p_1},\dots,\lambda_{p_{n-1-|J|}},\sum_{i\in
  I(\lambda)}\lambda_i)$ is sequential, and $mm((\lambda_{p_1},
\dots,\lambda_{p_{n-1-|J|}},\sum_{i\in I(\lambda)}\lambda_i))\geq 2$, 
we can conclude that the middle space in~\eqref{eq:4.13} is 
contractible, and hence $i$ in~\eqref{eq:4.12} is homotopic to 
a~trivial map. This yields the conclusion. 
\qed

\subsection{Resonances $(a^k,b^l,1^m)$.} $\,$

\begin{thm} \label{thm4.14}
  Let $a,b,k,l,m,r$ be positive integers, such that $b>1$, $m\geq r$,
  and $a=bl+r$.  Then
  \begin{equation}
    \label{eq:4.14}
    \cf(a^k,b^l,1^m)\simeq\susp(\cf(1^k)\wedge\cf(a,1^{m-r}))\vee
    (\cf(1^k)\wedge\cf(b^l,1^m)).
  \end{equation}
\end{thm}

\nin {\bf Note.} The restriction $m\geq r$ is unimportant. Indeed, if
$m<r$, then $a>bl+m$, hence $a$ is not involved in any resonance other
than $a=a$. This implies that
$\cf(a^k,b^l,1^m)=\cf(1^k)\times\cf(b^l,1^m)$, and we have determined
the homotopy type of $\cf(a^k,b^l,1^m)$ by the previous computations.

\vskip4pt

\nin {\bf Proof of the Theorem~\ref{thm4.14}.}

\nin Obviously, the condition $a>bl$ guarantees that the partition
$(a^k,b^l,1^m)$ is sequential, hence the Theorem~\ref{thm4.13} is
valid. It follows that if $k\geq 2$, then $\cf(a^k,b^l,1^m)$ is
contractible, hence~\eqref{eq:4.14} is true.

Furthermore, if $l\geq 2$, or, $l=1$ and $m\geq b$, then $(a,b^l,1^m)$
is strongly sequential, hence in this case~\eqref{eq:4.10} is valid,
which in new notations becomes precisely the equation~\eqref{eq:4.14}.

Finally, assume $l=1$ and $b>m\geq r\geq 1$. Let $a=b+d$. If
$\cf(a,1^{m-d})$ or $\cf(b,1^m)$ is contractible, then the map~$i$ in
the inclusion triple $\cf(a,1^{m-d})\stackrel{i}{\hra}\cf(b,1^m)\ra
\cf(b,1^m)/\cf(a,1^{m-d})$ is homotopic to a~trivial map, and we again
conclude~\eqref{eq:4.14}. If both of these spaces are not contractible
then $\cf(a,1^{m-d})\simeq S^{2y+\epsilon_2+1}$ and $\cf(b,1^m)\simeq
S^{2x+\epsilon_1+1}$, where nonnegative integers
$x,y,\epsilon_1,\epsilon_2$ are defined by
\begin{equation}
  \label{eq:4.15}
  m=bx+\epsilon_1,\,\,\,m-d=(b+d)y+\epsilon_2,\,\,\,
\epsilon_1,\epsilon_2\in\{0,1\}.
\end{equation}
Let us show that $2x+\epsilon_1>2y+\epsilon_2$. If $x>y$, then
$2x+\epsilon_1\geq 2x\geq 2y+2>2y+\epsilon_2$. From~\eqref{eq:4.15} we
have that $b(x-y)=d+dy+\epsilon_2-\epsilon_1$. If $x\leq y$, then the
left hand side is nonpositive. On the other hand, since $d\geq 1$, the
right hand side is nonnegative. Hence, both sides are equal to 0,
which implies $x=y$, $d=\epsilon_1=1$, $\epsilon_2=y=0$. This yields
$2x+\epsilon_1>2y+\epsilon_2$.

The homotopic triviality of the map $i$ follows now from the fact that
the homotopy groups of a~sphere are trivial up to the dimension of
that sphere, i.e., $\pi_k(S^n)=0$, for $0\leq k\leq n-1$.  
\qed

\subsection{Division chain resonances.} $\,$

We call the resonance $(b_n^{m_n},b_{n-1}^{m_{n-1}},\dots,b_1^{m_1})$
a~{\it division chain resonance} if $b_i\,\big|\,b_{i+1}$, for any $i\in
[n-1]$. For convenience, we assume $m_i\geq 1$, for $i\in [n]$, and
set $r_i=b_i/b_{i-1}$, for $n\geq i\geq 2$, and $r_1=b_1$.

Let us see that division chain resonances are strongly sequential.
First, we show that
$\lambda=(b_n^{m_n},b_{n-1}^{m_{n-1}},\dots,b_1^{m_1})$ is
sequential. Assume that
\begin{equation}\label{eq:dvr}
\sum_{i\in I}\alpha_i b_i=\sum_{j\in J}\beta_j b_j, 
\end{equation} 
and there are no equal size parts appearing on both sides. Set
$f=\max(I\cup J)$, $g=\min(I\cup J)$. We use induction on $f-g$. If
$f=g+1$ then the condition of sequentiality is obviously
satisfied. Otherwise, divide both sides by $b_g$. The number of parts
of size $1$ must be divisible by $r_{g+1}$, hence, in~\eqref{eq:dvr}
all the parts of size $b_g$ can be replaced by a~certain number of
parts of size $b_{g+1}$. By the induction assumption the condition of
sequentiality is satisfied for the new relation, hence it follows
for~\eqref{eq:dvr} as well.

Note that it also follows from the previous argument that $I(\lambda)$ 
must be of the form $\{p,p+1,\dots,n-mm(\lambda)-1,n-mm(\lambda)\}$,
for some~$p$.

It is now easy to see that $\lambda$ is strongly sequential. Assume
$b_n=b_{n-1}+\sum_{i\in I}\alpha_i b_i$, then
$(r_n-1)b_{n-1}=\sum_{i\in I}\alpha_i b_i$. The sequentiality
condition is true for the latter relation, hence the strong
sequentiality condition is true for the first one.

Thus, the Theorem~\ref{thm4.13} applies, and it yields:
\begin{enumerate}
\item if $m_n\geq 2$, then $\cf(\lambda)$ is contractible;
\item if $I(\lambda)$ exists, then
  \begin{multline}
    \label{eq:4.16}
  \cf(b_n,b_{n-1}^{m_{n-1}},\dots,b_1^{m_1}) \simeq
  (\cf(1)\wedge\cf(b_{n-1}^{m_{n-1}},\dots,b_1^{m_1}))\vee \\
  (S^1\wedge\cf(1)\wedge\cf(b_n,b_q^{\widetilde m_q},
  b_{q-1}^{m_{q-1}},\dots,b_1^{m_1})),
  \end{multline}
where $(b_n,b_q^{\widetilde m_q}, b_{q-1}^{m_{q-1}},\dots,b_1^{m_1})$
  is obtained from $(b_n,b_{n-1}^{m_{n-1}},\dots,b_1^{m_1})$ by
  removing the parts indexed by $I(\lambda)$. We have $\widetilde
  m_q\geq 1$.
\item If $I(\lambda)$ does not exist,
then
\begin{equation}
  \label{eq:4.17}
  \cf(b_n,b_{n-1}^{m_{n-1}},\dots,b_1^{m_1}) \simeq
  \cf(1)\wedge\cf(b_{n-1}^{m_{n-1}},\dots,b_1^{m_1}).
\end{equation}
\end{enumerate}

It is immediate from the formulae~\eqref{eq:4.16} and~\eqref{eq:4.17}
that each topological space
$\cf(b_n^{m_n},b_{n-1}^{m_{n-1}},\dots,b_1^{m_1})$ is homotopy
equivalent to a~wedge of spaces of the form $\cf(1)^\alpha\wedge
S^\beta$, where $\cf(1)^\alpha$ means an~$\alpha$-fold smash product
of $\cf(1)$. The natural combinatorial question which arises is how to
enumerate these spaces. We shall now construct a~combinatorial model:
a~weighted graph which yields such an~enumeration.

For convenience of notations, we set $m_0=1$. $\Gamma_\lambda$ is
a~directed weighted graph on the set of vertices $\{0,1,\dots,n\}$
whose edges and weights are defined by the following rule. For
$x,x+d\in\{0,\dots,n\}$, $d\geq 1$, there exists an~edge $e(x,x+d)$
(the edge is directed {\it from} $x$ {\it to} $x+d$) if and only if
$$b_{x+d}\,\big|\,b_{x+d-1}m_{x+d-1}+b_{x+d-2}m_{x+d-2}+\dots+
b_{x+1}m_{x+1}+b_x(m_x-1).$$
In this case the weight of the edge is defined as 
$$w(x,x+d)=(b_{x+d-1}m_{x+d-1}+\dots+
b_{x+1}m_{x+1}+b_x(m_x-1))/b_{x+d}.$$
Note that if $d\geq 2$ and there exists an edge $e(x,x+d)$, then there
exists an edge $e(x,x+d-1)$. 

We call a~directed path in $\Gamma_\lambda$ {\it complete} if it
starts in $0$ and ends in $n$. Let $\gamma$ be a~complete path in
$\Gamma_\lambda$ consisting of $t$ edges, $\gamma=(e(x_0,x_1),\dots,
e(x_{t-1},x_t))$, where $x_0=0$, and $x_t=n$. The weight of $\gamma$
is defined to be the pair $(l(\gamma),w(\gamma))$, where
$l(\gamma)=t$, and $w(\gamma)=\sum_{i=1}^t w(x_{i-1},x_i)$.

\begin{thm}\label{thm4.15}
  Let $\lambda=(b_n,b_{n-1}^{m_{n-1}},\dots,b_1^{m_1})$, then
  \begin{equation}
    \label{eq:4.18}
    \cf(\lambda)\simeq\bigvee_{\gamma}(\cf(1)^{l(\gamma)+w(\gamma)}
    \wedge S^{w(\gamma)}),
  \end{equation}
  where the wedge is taken over all complete paths of
  $\Gamma_\lambda$.
\end{thm}
\pr We use induction on $n$. The base of the induction is $n=1$.  In
 this case $\Gamma_\lambda$ is a~graph with only one edge $e(0,1)$,
 $w(0,1)=0$.  Thus, there is only one complete path. It has weight
 $(1,0)$, and $\cf(\lambda)\simeq\cf(1)$.

Next, we prove the induction step. We break up the proof in three
cases.

\vskip3pt

\nin {\bf Case 1.} {\it $I(\lambda)$ does not exist.}

\nin By~(\ref{eq:4.17}) we have
\begin{equation}
  \label{eq:4.19}
\cf(\lambda)\simeq\cf(1)\wedge\cf(b_{n-1}^{m_{n-1}},\dots,b_1^{m_1}).
\end{equation}
On the other hand, $I(\lambda)$ does not exist if and only if
$b_n>m_{n-1}b_{n-1}+\dots+m_1 b_1$. We also know that $n\geq 2$.  This
implies that there is at most one edge of the type $e(x,n)$, namely
$e(n-1,n)$. This edge exists if and only if $m_{n-1}=1$, in which case
$w(n-1,n)=0$.

If this edge does not exist then there are no complete paths in
$\Gamma_\lambda$ and, at the same time
$\cf(b_{n-1}^{m_{n-1}},\dots,b_1^{m_1})$ is contractible by the
previous observations. This agrees with~(\ref{eq:4.18}).

If, on the other hand, this edge does exist, then all complete paths
$\gamma$ must be of the type $\gamma=(\tilde\gamma,e(n-1,n))$, where
$\tilde\gamma$ is a~complete path from $0$ to $n-1$. Also in this 
case~(\ref{eq:4.19}) agrees with~(\ref{eq:4.18}).

\vskip3pt

\nin {\bf Case 2.} {\it $I(\lambda)$ exists and $m_{n-1}\geq 2$.}

\nin In this case $\cf(b_{n-1}^{m_{n-1}},\dots,b_1^{m_1})$ is
contractible, and 
\begin{equation}
  \label{eq:4.20}
  \cf(\lambda)\simeq S^1\wedge\cf(1)\wedge\cf(b_n,b_q^{\widetilde
  m_q}, b_{q-1}^{m_{q-1}},\dots,b_1^{m_1}),
\end{equation}
where $(b_n,b_q^{\widetilde m_q}, b_{q-1}^{m_{q-1}},\dots,b_1^{m_1})$
is as in~(\ref{eq:4.16}).

 Let $\tilde\lambda=(b_n,b_q^{\widetilde m_q},
b_{q-1}^{m_{q-1}},\dots,b_1^{m_1})$. We can describe the graph
$\Gamma_{\tilde\lambda}$: it is obtained from $\Gamma_\lambda$ by
\begin{enumerate}
\item removing all vertices indexed by $\{q+1,\dots,n-1\}$ and the
incident edges;
\item decreasing the weight of every existing edge $e(x,n)$ by $1$;
\item keeping all the existing edges with the old weights on the set
  $\{0,\dots,q-1,q\}$.
\end{enumerate}

This operation on $\Gamma_\lambda$ is well-defined, since there can be
no edges in $\Gamma_\Lambda$ of the type $e(x,n)$, for
$x\in\{q+1,\dots,n-1\}$, and since the weight of edges $e(x,n)$, for
$x\in\{0,\dots,q\}$ must be at least 1, as $\widetilde m_q\geq 1$.
Furthermore, it is clear from the above combinatorial description of
$\Gamma_{\tilde\lambda}$, that the set of the complete paths of
$\Gamma_{\tilde\lambda}$ is the same as that of $\Gamma_\lambda$, and
that the weights of the edges in these paths are also the same except
for the edge with the endpoint $n$, whose weight has been decreased
by~1. Thus,~(\ref{eq:4.20}) agrees with~(\ref{eq:4.18}) in this case.

\vskip3pt

\nin {\bf Case 3.} {\it $I(\lambda)$ exists and $m_{n-1}=1$.}

\nin This case is rather similar to the case 2, except that there is
an~edge $e(n-1,n)$ of weight $0$. Thus, $\Gamma_{\tilde\lambda}$
bookkeeps all the complete paths of $\Gamma_\lambda$, except for the
ones which have this edge $e(n-1,n)$.

However, the first term of the right hand side of~(\ref{eq:4.16})
bookkeeps the paths $(\tilde\gamma,e(n-1,n))$, just like in the
case~1.  Since the set of all complete paths of $\Gamma_\lambda$ is
the disjoint union of the sets of those paths which contain
$e(n-1,n)$, and those which do not, we again get that~(\ref{eq:4.16})
provides the inductive step for~(\ref{eq:4.18}).  \qed

\vskip4pt

\nin {\bf Examples.}

\nin (1) Let $\lambda=(a,1^l)$, for $a\geq 2$. Then $\Gamma_\lambda$
is a~graph on the vertex set $\{0,1,2\}$ having either one or two
edges:
\begin{enumerate}
\item it has in any case the edge $e(0,1)$, $w(0,1)=0$; 
\item if $a$ divides $l$, then it has the edge $e(0,2)$, in which case
  $w(0,2)=l/a$;
\item if $a$ divides $l-1$, then it has the edge $e(1,2)$, in which case
  $w(1,2)=(l-1)/a$.
\end{enumerate}
Clearly the Theorem~\ref{thm4.15} agrees with the
Theorem~\ref{thm4.2}. Indeed, if $\epsilon\not\in\{0,1\}$ (where
$\epsilon$ is taken from the formulation of the Theorem~\ref{thm4.2}),
then there are no complete paths in $\Gamma_\lambda$. If $\epsilon=0$,
then there is one path $(0,2)$ of weight $(1,l/a)$; and if
$\epsilon=1$, then there is one path $((0,1),(1,2))$ of weight
$(2,(l-1)/a)$. Thus,~(\ref{eq:4.18}) and~(\ref{eq:4.1}) are equivalent
in this case.

\nin (2) Let $\lambda=(8,4,2^3,1^6)$.
Then the graph $\Gamma_\lambda$ is

\[

  \begin{picture}(0,0)
    \includegraphics{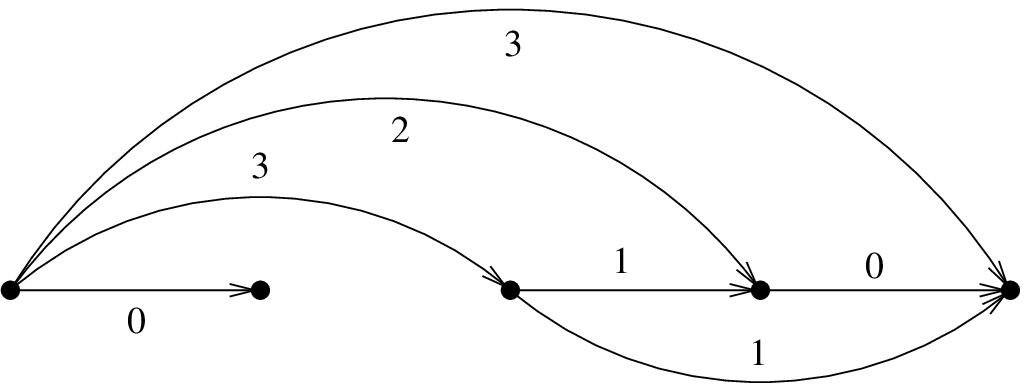}
  \end{picture}
  \begin{picture}(0,0)%
\includegraphics{graph.pstex}%
\end{picture}%
\setlength{\unitlength}{3947sp}%
\begingroup\makeatletter\ifx\SetFigFont\undefined%
\gdef\SetFigFont#1#2#3#4#5{%
  \reset@font\fontsize{#1}{#2pt}%
  \fontfamily{#3}\fontseries{#4}\fontshape{#5}%
  \selectfont}%
\fi\endgroup%
\begin{picture}(3714,1440)(994,-1181)
\end{picture}

  \]
\[
\text{Figure 1.}
\]

\noindent
It has 4 directed paths from $0$ to $4$ and, by the Theorem~~\ref{thm4.15}, 
we have 
$$\cf(\lambda)\simeq
(\cf(1)^3\wedge S^2)\vee
(\cf(1)^5\wedge S^3)\vee
(\cf(1)^6\wedge S^4)\vee
(\cf(1)^7\wedge S^4),$$ 
in particular $\Sigma_\lambda^{\mathbb R}\simeq S^5\vee
S^8\vee S^{10}\vee S^{11}$.

\section{Remarks on complexity of resonances.}

The main idea of all our previous computations was to find, for
a~given $n$-cut $S$, a~partition $\pi\in P(n)$, such that
$\lspan(S\setminus(\pi\dar S))\neq S$. Intuitively speaking, shrinking
the substratum corresponding to $\tilde\pi S$, where
$\un(\tilde\pi)=\pi$, {\it essentially} reduces the set of linear
identities in~$S$. It is easy to construct examples when such $\pi$
does not exist, e.q., Example~\ref{exsmres}(4).

These observations lead us to introduce a~formal notion of complexity
of a~resonance.

\begin{df} \label{df6.1}  $\,$

\nin
1) For $S\in\catr_n$, the {\bf complexity} of $S$ is denoted $c(S)$
and is defined by:
\begin{equation}
  \label{eq:6.1}
c(S)=\min\{|\Pi|\,|\,\Pi\subseteq P(n),\lspan(S\setminus(\Pi\dar S))\neq S\}.
\end{equation}

\nin
2) We define the complexity of an $n$-resonance to be the complexity
of one of its representing cuts. Clearly, it does not depend
on the choice of the representative.
\end{df}

\nin {\bf Note.} The number $c(S)$ would not change if we required the
partitions in $\Pi$ to have one block of size 2, and all other blocks 
of size~1.

\vskip4pt

The higher is the complexity of a resonance $[S]$, the less it is likely
that one can succeed with analyzing its topological structure using the
method of this paper. This is because one would need to take a~quotient
by a~union of $c([S])$ strata and it might be difficult to get a~hold
on the topology of that union.

\vskip3pt

We finish by constructing for an~arbitrary $n\in{\mathbb N}$,
a~resonance of complexity~$n$. Let
$\lambda_n=(a_1,\dots,a_n,b_1,\dots,b_n)$, such that
$a_i,b_i\in{\mathbb N}$, $a_i+b_j=a_j+b_i$, for $i,j\in[n]$, and all
other linear identities among $a_i$'s and $b_i$'s with coefficients
$\pm 1,0$ are generated by such identities. In other words, the cut
$S$ associated to $\lambda$ is equal to the set
\begin{equation}
  \label{eq:6.2}
  \Bigl\{(x_1,\dots,x_n,y_1,\dots,y_n)\in\{-1,0,1\}^{2n}\,\Big |\,
\sum_{i=1}^n y_i=0,\,\,x_i+y_i=0,\forall i\in[n]\Bigr\}.
\end{equation}

\nin
It is not difficult to construct such $\lambda_n$ directly:

\vskip3pt

\nin 1) Choose $a_1,\dots,a_n$, such that the only linear identities
with coefficients $\pm 1,0$ on the set $a_1,a_1,a_2,a_2,\dots,a_n,a_n$
are of the form $a_i=a_i$; in other words, there are no linear
identities with coefficients $\pm 2,\pm 1,0$ on the set
$a_1,\dots,a_n$. One example is provided by the choice $a_1=1$,
$a_2=3$, $\dots$, $a_n=3^{n-1}$.

\vskip3pt

\nin 2) Let $b_i=N+a_i$, for $i\in[n]$, where $N$ is sufficiently
large. As the proof of the Proposition~\ref{prop6.2} will show, it is
enough to choose $N>2\sum_{i=1}^n\lambda_i$. This bound is far from 
sharp, but it is sufficient for our purposes.

\begin{prop} \label{prop6.2}
  Let $S_n$ be the $n$-cut associated to the ordered sequence
of natural numbers $\lambda_n$ described above. Then $c(S_n)=n$.
\end{prop}

\pr First, let us verify that the cut $S_n$ associated to
$\lambda_n$ is equal to the one described in~\eqref{eq:6.2}. Take
$(x_1,\dots,x_n,y_1,\dots,y_n)\in S_n$. 

Assume first that $\sum_{i=1}^n y_i\neq 0$. Then,
$(x_1,\dots,x_n,y_1,\dots,y_n)$ stands for the identity
\begin{equation}
  \label{eq:6.3}
  \sum_{i\in I_1}a_i+\sum_{j\in J_1}b_j=
\sum_{i\in I_2}a_i+\sum_{j\in J_2}b_j,
\end{equation}
such that $|J_1|\geq|J_2|+1$. This implies that $N$ is equal to some 
linear combination of $a_i$'s with coefficients $\pm 2,\pm 1,0$.
This leads to contradiction, since $N>2\sum_{i=1}^n\lambda_i$.

Thus, we know that $\sum_{i=1}^n y_i=0$. Canceling $N\cdot|J_1|$ out
of~\eqref{eq:6.3} we get an~identity with coefficients $\pm 2,\pm 1,0$
on the set $a_1,\dots,a_n$. By the choice of $a_i$'s, this identity
must be trivial, which amounts exactly to saying that $x_i+y_i=0$, for
$i\in[n]$. 

Second, it is a~trivial observation that $c(S_n)\leq n$. Indeed, let
$\pi_i\in P(n)$ be a~partition with only one nonsingleton block
$(1,n+i)$, for $i\in[n]$. Then
$\lspan(S_n\setminus(\{\pi_1,\dots,\pi_n\}\dar S_n))\neq S_n$, since for any
$(x_1,\dots,x_n,y_1,\dots,y_n)\in S_n\setminus(\{\pi_1,\dots,\pi_n\}\dar S_n)$,
we have $x_1=0$.

Finally, let us see that $c(S_n)>n-1$. As we have remarked after the
Definition~\ref{df6.1}, it is enough to consider the case when the
partitions of $\Pi$ have one block of size 2, and the rest are
singletons. Let us call the identity $a_i+b_j=a_j+b_i$ {\it the 
elementary identity indexed $(i,j)$}. 

From the definition of the closure operation $\dar$ it is clear that
an~elementary identity indexed $(i,j)$ is not in $S_n\setminus(\Pi\dar
S_n)$ if and only if the partition whose only nonsingleton block is
$(i,n+j)$ belongs to $\Pi$, or the partition whose only nonsingleton
block is $(j,n+i)$ belongs to $\Pi$. That is because the only reason
this identity would not be in $S_n\setminus(\Pi\dar S_n)$ would be
that one of these two partitions is in $\Pi\dar S_n$. But, if such
a~partition is in $\Pi\dar S_n$, then it must be in $\Pi$: moves (2)
of the Definition~\ref{df3.1} can never produce a~partition whose only
nonsingleton block has size~2, while the moves (3) of the
Definition~\ref{df3.1} may only interchange between partitions
$(i,n+j)$ and $(j,n+i)$ in our specific situation. Thus, we can
conclude that if $|\Pi|\leq n-1$, then at most $n-1$ elementary
identities are not in $S_n\setminus(\Pi\dar S_n)$.

Next, we note that for any distinct $i,j,k\in[n]$, the elementary
identities $(i,j)$ and $(j,k)$ imply the elementary identity $(i,k)$.
Let us now think of elementary identities as edges in a~complete graph
on $n$~vertices, $K_n$. Then, any set $M$ of elementary identities
corresponds to a~graph $G$ on $n$ vertices, and the collection of the
elementary identities which lie in the $\lspan M$ is encoded by the
{\it transitive closure} of $G$. It is a~well known combinatorial fact
that $K_n$ is $(n-1)$-connected, which means that removal of at most
$n-1$ edges from it leaves a~connected graph. Hence, if we remove at
most $n-1$ edges from $K_n$ and then take the transitive closure, we
get $K_n$ again. Thus, if $|\Pi|\leq n-1$, all elementary identities
lie in $\lspan(S_n\setminus(\Pi\dar S_n))$. Since the elementary
identities generate the whole $S_n$, we conclude that
$S_n=\lspan(S_n\setminus(\Pi\dar S_n))$, hence $c(S_n)>n-1$.
\qed

\end{document}